\documentclass[11pt]{article}
\usepackage{amssymb}
\usepackage{amsmath}
\usepackage{amsthm}

\makeatletter
\def\LaTeX{\leavevmode L\raise.42ex
    \hbox{\kern-.3em\size{\sf@size}{0pt}\selectfont A}\kern-.15em\TeX}
\makeatother

\newcommand{\BibTeX}{{\rm B\kern-.05em{\sc
          i\kern-.025emb}\kern-.08em\TeX}}

\makeatletter
\def\@currentlabel{2.1}\label{e:dispaa}
\def\@currentlabel{2.21}\label{e:dispau}
\def\@currentlabel{2.22}\label{e:dispav}
\def\@currentlabel{2.23}\label{e:dispaw}
\def\@currentlabel{2.24}\label{e:dispax}
\def\theequation{\thesection.\@arabic\c@equation}
\makeatother

\newcounter{mnotecount}[section]

\newcommand{\rmnote}[1]{}

\renewcommand{\theequation}{\arabic{section}.\arabic{equation}}
\newtheorem{thm}{Theorem}[section]
\newtheorem{lem}[thm]{Lemma}
\newtheorem{cor}[thm]{Corollary}
\newtheorem{prop}[thm]{Proposition}
\theoremstyle{definition}

\newtheorem{rem}[thm]{Remark}

\newcommand{\B}{\mathbb B}

\newcommand{\R}{\mathbb{R}}
\newcommand{\N}{\mathbb{N}}

\newcommand{\aint}{-\!\!\!\!\!\!\int}

\def \p{\partial}

\begin{document}
\title{Remarks on two fourth order elliptic problems in whole space}
\author{Baishun Lai and Dong Ye}
\date{}
\maketitle
\begin{abstract}
We are interested in entire solutions for the
semilinear biharmonic equation $\Delta^{2}u=f(u)$ in $\R^N$, where $f(u)=e^{u}$ or $-u^{-p}\ (p>0)$. For the exponential case,
we prove that for the polyharmonic problem $\Delta^{2m}u = e^u$ with positive integer $m$, any classical entire solution verifies $\Delta^{2m - 1} u< 0$, this completes the results in \cite{Dupaigne, xu-wei};
we obtain also a refined asymptotic expansion of radial separatrix solution to $\Delta^{2}u= e^u$ in $\R^3$, which answers a question in \cite{Berchio}. For the negative power case,
we show the nonexistence of the classical entire solution for any $0<p\leq 1$.
\end{abstract}

\noindent
{\small {\bf Mathematics Subject Classification (2000):} 35J91, 35B08, 35B53, 35B40.

}

\setcounter{equation}{0}
 \setcounter{equation}{0}
\section{Introduction}
In the present note, we are interested in entire solutions for two
semilinear biharmonic equations
\begin{equation}\label{E1.1}
\Delta^{2}u=e^{u} \ \ \mbox{in}\ \R^{N}
\end{equation}
and
\begin{align}
\label{bi}
\Delta^2 u = -u^{-p}\;\; \mbox{in }\R^N, \quad \mbox{where $p > 0$}.
\end{align}
Recently, the fourth order equations have attracted the interest of many researchers. In
particular, a lot of efforts have been devoted to understand the existence, multiplicity, stability and qualitative properties of solutions for $\Delta^2u = f(u)$ with classical nonlinearities, like the polynomial growth $f(u) = u^p$, the exponential growth $f(u) = e^u$ and the negative power situation $f(u) = -u^{-p}$. For equation \eqref{E1.1}, in the conformal dimension $N=4$, \eqref{E1.1} appears naturally in conformal geometry as the constant $Q$-curvature problem, the
existence and asymptotic behaviour of solutions with finite total curvature, i.e.~$e^u \in L^1(\R^4)$ were studied in \cite{Chang,Lin,Wei}. Entire radial solutions of \eqref{E1.1} were also studied for $N\geq5$ in \cite{Arioli} and the stability of these entire radial solutions were
considered in \cite{Berchio, Dupaigne}. In particular, it is proved by \cite{Berchio} that
\eqref{E1.1} admits no radial entire solution if $N=2$.

\medskip
Recently, Farina informed us that a very general nonexistence result was proved by Walter in 1957, see \cite{Wa}. In particular, Walter proved that no classical entire solution exists in $\R^2$ for the polyharmonic problem $\Delta^{2m}u = e^u$ with any positive integer $m$.  Here we give an alternative proof (see Corollary \ref{C1.2} and Remark \ref{new1} below). Indeed, we will make use of a general observation for entire solutions to $\Delta^{2m}u = e^u$. By classical or smooth solution to $\Delta^\ell u = f(u)$ with $\ell \in \N^*$, we mean a solution in the class $C^{2\ell}$, equivalently all $2\ell$-th order derivatives of $u$ are continuous.

\begin{thm}\label{Th1.1}
Let $u$ be a classical solution of $\Delta^{2m}u = e^u$ in $\R^N$ with $m \in \N^*$, then $\Delta^{2m-1} u < 0$, i.e.~$(-\Delta)^{2m-1} u > 0$ in $\R^N$.
\end{thm}
We note that similar results were obtained by \cite{Dupaigne, xu-wei} under additional conditions.
The authors in \cite{Dupaigne} considered solutions to \eqref{E1.1} which are stable outside a bounded domain. In \cite{xu-wei}, it was proved that $(-\Delta)^{\ell-1} u > 0$ for any classical entire solution of $(-\Delta)^\ell u = e^u$ with $\ell \geq 2$, satisfying $u(x)=o(|x|^{2})$ at infinity.

\medskip
It is worthy to mention that the corresponding result is no longer true for classical entire solutions to $(-\Delta)^\ell u = e^u$ with odd $\ell$. In fact, Farina and Ferrero prove that for any $m \geq 1$, there are infinitely many entire radial solutions of $(-\Delta)^{2m+1} u = e^u$ such that $\Delta^{2m} u$ changes sign, see Lemma 6.8 and the proof of Lemma 5.4 in \cite{FF}. See also \cite{Wa1} for entire radial solutions of the equation $\Delta^\ell u = e^u$ with $\ell > 1$, $N  \geq 3$.

\medskip
On the other hand, for $N\geq3$, it is known that \eqref{E1.1} admits infinitely many smooth radial solutions. These radial solutions are
of either exactly quadratic growth or logarithmic growth at infinity
for $N \geq 4$ (see \cite{Arioli, Berchio}). For $N=3$, it is proved in \cite{Berchio} that the radial solution is
of either exactly quadratic growth or it verifies $u(r)\leq -Cr$ at
infinity for some $C>0$.  More precisely, let $u_{\alpha, \beta}$ be the unique radial solution of
\begin{align}\label{E1.2}
\left \{
\begin{array}{ll}
 \Delta^{2} u_{\alpha,\beta}(r)= e^{u_{\alpha,\beta}(r)} \ \mbox{for}\ \ r\in [0, R(\alpha,\beta)),\\
 u_{\alpha,\beta}(0)=\alpha,\; \Delta u_{\alpha,\beta}(0)=\beta, \; u'_{\alpha,\beta}(0)=(\Delta u_{\alpha,\beta})'(0)=0,
 \end{array}
\right.
\end{align}
where $[0,R(\alpha,\beta))$ denotes the maximal interval of existence. Noting that the equation \eqref{E1.2}
is invariant under the scaling transformation
$$
u_{\lambda}(x)=u(\lambda x)+4\ln \lambda,  \;\lambda>0.
$$
Therefore, we need only to understand the case $\alpha=0$. We will denote $u_{0,\beta}$ by $u_{\beta}$ and
$R(0,\beta)$ by $R(\beta)$ for simplicity. It has been proved in \cite{Arioli,Berchio} that any local solutions to \eqref{E1.2} satisfies
\begin{equation}\label{E1.3}
u_{\beta}(r)\geq \frac{\beta}{2N}r^{2}\ \ \mbox{for all }\ r\in [0, R(\beta)).
\end{equation}
Furthermore, there exists $\beta_0 \in (-\infty, 0)$ such
that
\begin{enumerate}
\item[(i)] For $\beta<\beta_0$, then $R(\beta)=+\infty$ and in
addition to \eqref{E1.3}, one has the upper bound
$$
u_{\beta}(r)\leq -\frac{\beta_{0}-\beta}{2N}r^{2}\ \ \mbox{for all }\ r\in [0, \infty);
$$
\item[(ii)] For $\beta=\beta_0$, the solution $u_{\beta_0}$, called separatrix verifies
$$
\left \{
\begin{array}{ll}
   u_{\beta_{0}}(r)\leq -Cr  ,&\ \mbox{if}\  N=3  \mbox{ and $r$ large, with } C>0;\\
  u_{\beta_{0}}(r) =-4\ln\left(1+\frac{e^{\frac{\alpha}{2}}}{8\sqrt{6}}r^{2}\right),& \  \mbox{for}\ N=4;\\
  \lim_{r\to \infty}\left[u_{\beta_0}(r)+4\ln r\right]=\ln[8(N-2)(N-4)],& \  \mbox{for}\ N\geq 5.
  \end{array}
\right.
$$
\item[(iii)] For $\beta>\beta_{0}$, $R(\beta)<\infty$ and $\lim_{r\nearrow R(\beta)}u_{\beta}(r)=\infty$.
\end{enumerate}

An open problem was left for the exact asymptotic behaviour of the separatrix $u_{\beta_0}$
in dimension three, see \cite{Berchio}. The following result answers this issue.
\begin{thm}\label{Th1.3}
Let $\beta_0$ be defined as above and $N=3$. Then we have, as $r\to\infty$, $u_{\beta_0}(r)= \alpha_{1}r + \alpha_2 + \alpha_3r^{-1} + O(e^{-cr})$ where $c > 0$ and
\begin{align*}
\alpha_{1} = \frac{-1}{8\pi}\int_{\R^3}e^{u_{\beta_0}} dx, \;\; \alpha_2 = \frac{1}{8\pi}\int_{\R^3}|x|e^{u_{\beta_0}} dx, \;\; \alpha_3 = \frac{-1}{24\pi}\int_{\R^3}|x|^2e^{u_{\beta_0}} dx.
\end{align*}
\end{thm}

The second part of the note is devoted to consider the classical solutions of equation \eqref{bi}. Recently, the radial solutions to \eqref{bi} are studied in \cite{DFG}, and some Liouville type results are obtained for stable entire solutions of \eqref{bi} in \cite{GW}. We can remark that all these results concern the negative exponent $-p$ with $p > 1$, and it seems curious for us that no study existed for entire solutions of \eqref{bi} with $p \leq 1$. Here we prove that no such entire solution could exist if $p \in (0, 1]$, that is
\begin{thm}\label{Th1.5}
If $0 <p\leq 1$, the equation \eqref{bi} admits no entire smooth solution.
\end{thm}
In fact, our proof is inspired by the work of Choi-Xu in \cite{cx}, where the above result has been established for $N = 3$.

\setcounter{equation}{0}
 \setcounter{equation}{0}
\section{Proof of Theorem \ref{Th1.1}}
In this section,  we prove Theorem \ref{Th1.1}. In the following, for a given function $f$, we write $$\overline{f}(r)=\aint_{\partial B_{r}(0)}f d\sigma = \frac{1}{|\p B(0, r)|}\int_{\partial B_{r}(0)}f d\sigma, \quad  \forall\; \; r > 0,$$
where $|\p B(0, r)|$ denotes the volume of the sphere. Furthermore, we will consider $\Delta^{2m} u = e^u$ as a
system:
\begin{align}
\label{defv}
v_1 := u, \;\; v_{k+1} := \Delta v_k \; \mbox{for $1 \leq k \leq 2m - 1$ so that } \Delta v_{2m}=e^{u}  \ \ \mbox{in}\ \R^{N}.
\end{align}

\medskip
\noindent
{\bf Proof of Theorem \ref{Th1.1}.} First we show that $v_{2m} = \Delta^{2m-1} u \leq 0$. If it is not the case, there is
a point $x_0 \in R^{N}$ such that $v_{2m}(x_0) > 0$. Up to a translation, we may assume that $x_0 =0$. Therefore with $v_k$ given by \eqref{defv}, $\overline{v_k}(r)$ satisfy
\begin{align}\label{E2.2}
\Delta \overline{v_k}= \overline{v_{k+1}} \; \mbox{for $1 \leq k \leq 2m - 1$}, \quad
\Delta \overline{v_{2m}}=\overline{e^{u}}\geq e^{\overline{u}}  \ \ \mbox{in}\ \R^{N}.
\end{align}
Remark that $\Delta\overline{v_{2m}}= r^{1-N}(r^{N-1}\overline{v_{2m}}')' =\overline{e^{u}} >0$, so $\overline{v_{2m}}$ is increasing w.r.t.~the radius $r$. There holds $\Delta
\overline{v_{2m-1}}\geq \overline{v_{2m}}(0) >0$. Integrating it, we get
\begin{equation*}
v_{2m-1}(r)\geq v_{2m-1}(0)+ \frac{\overline{v_{2m}}(0)}{2N}r^{2}.
\end{equation*}
Hence $v_{2m-1}(r)\to \infty$ as $r\to\infty$. By iteration, we see that $\overline u(r) = \overline v_1(r)\to \infty$ as $r\to\infty$. Now Let $r=e^{t}, w(t)=\overline{u}(e^{t})$, direct calculation yields
\begin{equation}
\label{eq2.4}
e^{4mt} e^{w(t)} = e^{4mt}e^{\overline u(r)} \leq e^{4mt}\Delta^{2m}\overline u(r) = w^{(4m)}(t) + \sum_{i = 1}^{4m-1} c_i w^{(i)}(t)
\end{equation}
where $c_i$ are some constants depending only on $N$ and $i$. Here and after, $g^{(i)}$ denotes the $i$-th derivative of a function $g$. Since $\lim_{t\to\infty}w(t)=\infty$, there exists $T_1$ such that
$$
e^{4mt} e^{w(t)}\geq w^{2}(t)\; \mbox{ for all }\; t \geq T_1.
$$
We apply now the test function method developed by  Mitidieri and Pohozaev in \cite{Mit}. More precisely, we can choose a nonnegative function $\phi_0\in C_{0}^{\infty}[0,\infty)$ satisfying $\phi_0 > 0$ in $[0, 2)$,
$$
\phi_{0}(\tau)=\left \{
\begin{array}{ll}
1\ & \ \mbox{for}\ \tau\in [0, 1]\\
0 & \ \mbox{for}\ \tau \geq2.  \\
\end{array}
\right.
\quad  \mbox{and}\quad
\int_{0}^{2}\frac{|\phi_{0}^{(i)}(\tau)|^{2}}{\phi_{0}(\tau)}d\tau:=A_{i}<\infty\;\; \forall\; i\in \N.
$$
Let $T > T_1$, multiplying \eqref{eq2.4} by $\phi(t)=\phi_{0}\left(\frac{t-T_{1}}{T-T_{1}}\right)$ and integrating by parts, we obtain
\begin{align}\label{E2.6}
\int_{T_{1}}^{\infty}\left[\phi^{(4m)}(t)+ \sum_{i = 1}^{4m-1} (-1)^ic_i\phi^{(i)}(t)\right]w(t)dt \geq \int_{T_{1}}^{\infty}w^{2}(t)\phi(t) dt - C.
\end{align}
By Young's inequality, for any $\epsilon>0$,  $ \exists \ C_\epsilon>0$ such that
$$
w(t)\phi^{(i)}(t)\leq \epsilon w^{2}(t)\phi(t)+C_\epsilon\frac{|\phi^{(i)}(t)|^{2}}{\phi(t)}, \quad \forall\; t \in [T_1, 2T-T_1).
$$
Then, provided that $\epsilon$ is chosen sufficiently small, \eqref{E2.6}
yields
\begin{align*}
C'\sum_{i=1}^{4m}A_i(T-T_{1})^{1-2i}=C'\sum_{i=1}^{i=4m}\int_{T_{1}}^{2T-T_1}\frac{|\phi^{(i)}(t)|^{2}}{\phi(t)}dt& \geq
\int_{T_{1}}^{2T-T_1}w^{2}(t)\phi(t)dt - C'',\\
& \geq \int_{T_{1}}^{T}w^{2}(t)dt - C'',
\end{align*}
with fixed constants $C', C'' > 0$. Let $T\to\infty$, we observe a contradiction with $w(t)\to
\infty$. So we have $v_{2m} \leq 0$ in $\R^N$.

\medskip
Now suppose that there exists $x_0 \in \R^N$ verifying $v_{2m}(x_0) = 0$, then $x_0$ is a maximum of $v_{2m}$, hence
$\Delta v_{2m}(x_0) \leq 0$ which is just impossible as $\Delta v_{2m} = e^u$, so $\Delta^{2m - 1} u = v_{2m} < 0$ in $\R^N$.\qed

\medskip
As an immediate consequence of Theorem \ref{Th1.1}, we can claim
\begin{cor}\label{C1.2}
For any $m \in \N^*$, the equation $\Delta^{2m} u = e^u$ admits no classical entire solution in $\R^2$.
\end{cor}

\noindent
{\bf Proof.} We suppose by contradiction that $u$ is
a smooth function verifying $\Delta^{2m} u = e^u$ in $\R^2$. Using Theorem \ref{Th1.1}, $v := \Delta^{2m - 1}u < 0$ in $\R^2$. Moreover,
$$
\overline v'(r)=\frac{1}{2\pi r}\int_{\B_{r}(0)}\Delta \overline v dx = \frac{1}{2\pi r}\int_{\B_{r}(0)}\Delta^{2m} u dx = \frac{1}{2\pi r}\int_{\B_{r}(0)} e^u dx\geq \frac{C}{r}, \quad \forall\; r\geq 1,
$$
where $C$ is a positive constant. Hence
$$
\overline{v}(r)-\overline{v}(1)=\int_1^{r}\overline{v}'(r)dr\geq C\ln r, \quad \forall\; r\geq 1.
$$
This contradicts the fact $\overline v(r)<0$ if we tend $r$ to $\infty$, so we are done. \qed

\begin{rem}
\label{new1}
By adapting similar approach,
the results of Theorem \ref{Th1.1} and Corollary  \ref{C1.2} hold true for the equation $\Delta^{2m} u = f(u)$ with
general convex, positive nonlinearity $f$ verifying
\begin{align}
\label{growth}
\liminf_{t\to\infty}
f(t)t^{-1-\mu} > 0 \;\;\mbox{
 for some $\mu > 0$}.\end{align}
We should mention that Walter proved in \cite{Wa} the nonexistence of smooth entire solution to $\Delta^{2m} u = f(u)$ in $\R^2$ for any $m \in \N^*$ and any positive function $f$ satisfying \eqref{growth}, without the convexity assumption.
\end{rem}

 \setcounter{equation}{0}
\section{Proof of Theorem \ref{Th1.3}}
We will use here the notations in Introduction for radial solutions, and also the results (i)-(iii)
cited there, given by \cite{Arioli, Berchio}. Recall that $u_\beta$ is the unique radial solution of
\begin{align}
\label{2.1n}
\Delta^2 u_\beta = e^{u_\beta}, \;\; \Delta u_\beta(0)=\beta, \;\; u_\beta(0) = u'_\beta(0)=(\Delta u_\beta)'(0)=0;
\end{align}
and the solution exists globally if and only if $\beta \leq \beta_0$. First, we show the following
characterization of the separatrix solution $u_{\beta_0}$.
\begin{prop}
\label{P3.1}
For any $\beta \leq \beta_0$, $\lim_{r\to\infty}\Delta u_{\beta}(r)\leq0$ and
$\lim_{r\to\infty}\Delta u_{\beta}(r)=0$ if and only if
$\beta=\beta_{0}$.
\end{prop}

\noindent
{\bf Proof.} For any solution $u$ of \eqref{E1.1}, $$\frac{d \Delta u(r)}{dr}=r^{N-1}\int_{0}^{r}s^{1-N}e^{u}ds>0.$$
According to Theorem \ref{Th1.1}, $\lim_{r\to\infty}\Delta
u_{\beta}(r)= \sigma\leq 0$ exists. For $\beta<\beta_{0}$, we see that $\sigma<0$, since $u_\beta \leq -Cr^2$ by (i)
and $\sigma=0$ implies readily that $u_{\beta}(r)=o(r^{2})$ at $\infty$.

\medskip
Similarly, we easily obtain $\lim_{r\to\infty}\Delta u_{\beta_{0}}=0$ for $N \geq 4$ by (ii).
Consider now $u_{\beta_0}$ when $N=3$. In fact, we will prove that if $\sigma < 0$, then $\beta < \beta_0$.

\medskip
For $N = 3$, \eqref{E1.2} reads
\begin{equation}\label{E3.1}
(r^{4}u'''(r))'=r^{4}e^{u},\quad \forall \; r>0.
\end{equation}
Integrating over $[0, r]$, we see that for all $r\geq 1$,
$$
r^{4}u'''(r) =\int_{0}^{r}s^{4}e^{u(s)}ds \leq \int_{0}^{\infty}s^{4}e^{u(s)}ds < \infty.
$$
Here we used the fact that $u(r)\leq-Cr$ for $r$ large. Thus $u'''(r)< Cr^{-4}$ for $r\geq 1$. Suppose now $\sigma = \lim_{r\to\infty}\Delta
u(r) <0$ for some entire solution $u$ of \eqref{2.1n} with $N= 3$.
As $$u'(r)=r^{-2}\int_{0}^{r}s^{2}\Delta u(s)ds,$$ we have then
$$
u(r)\sim \frac{\sigma}{6}r^{2},\;\; u'(r)\sim \frac{\sigma}{3}r,\;\; u''(r)\sim \frac{\sigma}{3} \quad \mbox{when $r \to \infty$}.
$$
Consider now the function $\tilde{u}$ defined by
$$
\tilde{u}(r)=-\epsilon r^{2}+\ln(1+r)-b
$$
where
$$\epsilon > 0, \quad b\geq \ln\left(\max_{\R_+}\psi\right) \;\; \mbox{with }\; \psi(r):=\frac{r(1+r)^{5}}{2(r+4)}e^{-\epsilon
r^{2}} \;\mbox{in }  \R_+.$$
Direct
computation shows that $\tilde{u}$ is supersolution of
\eqref{E3.1} in $\R^3$ and
$$\tilde{u}'(r)=-2\epsilon r+\frac{1}{r+1},\;\;
\tilde{u}''(r)=-2\epsilon-\frac{1}{(r+1)^{2}},\;\;
\tilde{u}'''(r)=\frac{2}{(r+1)^{3}}.$$ Hence, if we fix $\epsilon \in (0, -\sigma/6)$ and some large enough
$r_{0}$, there hold $u^{(i)}(r_{0})<\tilde{u}^{(i)}(r_{0})$ for $0 \leq i \leq 3$. By continuous dependence on initial data, there is $\beta_{1}>\beta = -\Delta u(0)$ such that $
u_{\beta_1}^{(i)}(r_{0})<\tilde{u}^{(i)}(r_{0})$ for $0\leq i \leq 3$. We claim then
\begin{align}\label{2.2n}
u_{\beta_1}(r)< \tilde{u}(r) \quad \mbox{for all }  r \geq r_0.
\end{align}
If it is not the case, then
$$r_1 = \sup\left\{s > r_0 \mbox{  s.t. } u_{\beta_1}(r)<
\tilde{u}(r) \mbox{ in } [r_0, s]\right\} < \infty.$$
By \eqref{E3.1}, we have $(r^{4}u_{\beta_{1}}'''(r))'< (r^{4}\tilde{u}'''(r))'$ in $[r_{0}, r_1)$, and successive integrations yield that $u_{\beta_1}'<\tilde{u}'$ on
$[r_{0}, r_1)$, hence $u_{\beta_1}(r_1)<\tilde{u}(r_1)$ . This contradicts the definition of $r_1$, so the claim \eqref{2.2n} holds true. By the point (iii), $u_{\beta_{1}}$ is defined then for all $r\geq 0$ which means that $\beta_{1}\leq \beta_{0}$, so $\beta<\beta_{0}$. \qed

\medskip
\noindent
{\bf Proof of Theorem \ref{Th1.3}.} To simplify the presentation, we erase the index
$\beta_0$ and denote $u_{\beta_0}$ by $u$. Recall that $u \leq -Cr$ for some $C > 0$ by (ii). Let $v= -\Delta u$, then we have
$$
v(r)=\beta_0 - \int_{0}^{r}s^{-2}\int_{0}^{s}t^{2}e^{u(t)}dtds, \quad \forall \; r > 0.
$$
Applying Proposition \ref{P3.1}, as $\lim_{r \to \infty} v(r) = 0$, we get
\begin{align*}
\begin{split}
v(r)= \int_{r}^{\infty}s^{-2}\int_{0}^{s}t^{2}e^{u(t)}dtds & = \frac{1}{r} \int_{0}^{r}t^{2}e^{u(t)}dt + \int_{r}^{\infty}te^{u(t)}dt\\
&= \frac{1}{4\pi r}\int_{\R^{3}}e^{u}dx -\frac{1}{r}\int_{r}^{\infty}t^{2}e^{u}dt + \int_{r}^{\infty}te^{u}dt.
\end{split}
\end{align*}
Therefore
\begin{equation}\label{E3.2}
(r^{2}u'(r))'=ar+r\int_{r}^{\infty}t^{2}e^{u}dt-r^{2}\int_{r}^{\infty}te^{u}dt \;\; \mbox{where } \; a=-\frac{1}{4\pi}\int_{\R^3}e^{u}dx.
\end{equation}
Integrating \eqref{E3.2}, we obtain
\begin{align*}
\begin{split}
u(r)=\frac{ar}{2}+\frac{1}{2}\int_{0}^{r}t^{3}e^{u}dt-\frac{1}{6r}\int_{0}^{r}t^{4}e^{u}dt
+\frac{r^{2}}{6}\int_{r}^{\infty}te^{u}dt.
\end{split}
\end{align*}
Then it is easy to get the claimed expansion for $u$. \qed

\setcounter{equation}{0}
\section{Proof of Theorem \ref{Th1.5}}
The proof of Theorem \ref{Th1.5} is based on the following lemma.
\begin{lem}
\label{l1.3}
If $u$ is a smooth solution of \eqref{bi}, then $\Delta u>0$ in
$\R^{N}$.
\end{lem}

 Indeed, this Lemma is an immediate consequence of the followin result.
\begin{lem}\label{l1new}
If $u$ is a $C^4$ lower bounded function verifying that $\Delta^2 u < 0$ in $\R^N$, then $\Delta u>0$ in
$\R^{N}$.
\end{lem}

\noindent
{\bf Proof.} First we show by contradiction that $\Delta u\geq 0$. Suppose that there is  $x_0\in \R^N$ verifying
$\Delta u(x_0)<0$. By translation, we can assume that $x_0 = 0$. Let $w = \Delta u$, then $\Delta \overline{u}=\overline{w}$ and $\Delta\overline{w} = \overline{\Delta^2 u} < 0$ where $\overline u$ and $\overline w$ are the average over sphere for $u$ and $w$. Consequently $\overline{w}'(r)\leq 0$, hence $\overline{w}(r)\leq \overline{w}(0) = \Delta u(0) <0$. Therefore
$\Delta \overline{u}\leq \overline{w}(0)$ in $\R^N$ which yields
$$
\overline{u}(r)\leq \overline{u}(0)+\frac{\overline{w}(0)}{2N}r^{2}
$$
We get $\overline{u}(r)<0$ for $r$ large enough, which is impossible since $u$ is lower bounded. So $\Delta u\geq 0$ in $\R^N$. Now if there is $x_1\in \R^{N}$ such that $\Delta u(x_1)=0$. Thus $x_1$ is a minimum point of $\Delta
u$ and $\Delta^2 u(x_1)\geq 0$, which contradicts the hypothesis, so the proof is completed. \qed

\medskip
From the above proof, as $\overline w \leq w(0)$, we immediately have
\begin{cor}\label{C1.4}
If $u$ is a $C^4$ lower bounded solution in $\R^N$ verifying $\Delta^2 u < 0$ in $\R^N$, then there exists $C > 0$ such that $\overline{u}(r)\leq C(1 + r^2)$ for any $r \geq 0$.
\end{cor}

\noindent
{\bf Proof of Theorem \ref{Th1.5}.} For $N=1$, we have $u''>0$ from
Lemma \ref{l1.3} and $u^{(4)}<0$. However, except being constant, any function cannot be
concave and lower bounded on $\R$, so we get the nonexistence of entire solution for $u^{(4)} = u^{-p}$ in $\R$ for any $p > 0$. For $N=2$, the superharmonic function $\Delta u$ is bounded from below by Lemma \ref{l1.3}, so it must be constant,  again it cannot verify the \eqref{bi}, so we are done.

\medskip
Consider from now on $N\geq 3$,
we claim that if $u$ is a smooth solution of \eqref{bi},
then
\begin{equation}\label{Eq1.3}
\mbox{there exists $C > 0$ such that } \overline u(r)\geq C r^{\frac{4}{p+1}}, \; \forall \; r > 0.
\end{equation}
In fact, $\overline w$ is decreasing where $w = \Delta u$, and $\overline u$ is increasing as $\overline w > 0$ by Lemma \ref{l1.3}. Using $\Delta \overline u = \overline w$, we have, by the monotonicity of $\overline w$,
\begin{equation}\label{Eq1.4}
\overline{u}(r)\geq u(0)+\frac{\overline{w}(r)}{2N}r^{2}.
\end{equation}
On the other hand, By Jensen's inequality,
$$
f(r) := -\Delta {\overline w} (r) = \overline{u^{-p}}(r) \geq {\overline u}^{-p}(r) > 0.
$$
For any $s \geq r > 0$,
\begin{align*}
{\overline w}'(s) = -s^{1-N}\int_0^s t^{N-1}f(t) dt \leq -s^{1-N}\int_0^r t^{N-1}f(t) dt,
\end{align*}
so we get, using the monotonicity of $\overline u$,
\begin{align}
\label{1.2new}
\begin{split}
{\overline w}(r) \geq {\overline w}(2r) + \int_r^{2r} s^{1-N}\int_0^r t^{N-1}f(t) dt ds & \geq {\overline w}(2r) + Cr^{2-N}\int_0^r t^{N-1}f(t) dt\\
& \geq Cr^{2-N}\int_0^r t^{N-1}{\overline u}^{-p}(t) dt\\
& \geq Cr^2{\overline u}^{-p}(r),
\end{split}
\end{align}
Inserting into \eqref{Eq1.4}, we have
$$\overline u(r) \geq u(0) + Cr^4{\overline u}^{-p}(r) \geq  C r^4{\overline u}^{-p}(r).$$
Hence (\ref{Eq1.3}) follows.

\medskip
Combining (\ref{Eq1.3}) and Corollary \ref{C1.4}, if $u$ is a classical solution of \eqref{bi}, necessarily there holds $p\geq1$.
Finally, we will exclude the case $p=1$. Let $u$ be a smooth entire solution to $\Delta^2 u = -u^{-1}$, then $\overline{u}$ is a subsolution to the following equation
\begin{equation}\label{Eq1.5}
 \Delta^{2}U(r)+U^{-1}(r)=0, \;\; U(0)=u(0),\; U''(0)= \overline{u}''(0),\; U'(0) = U'''(0)=0.
\end{equation}
Consider
$$
Z(r) = u(0)+\frac{\overline{u}''(0)}{2}r^{2}.
$$
Obviously, $Z$ is biharmonic and a supersolution of (\ref{Eq1.5}). A comparison
principle (see Lemma 3.2 in \cite{Mc}) ensures that $Z \geq \overline{u}$, and there is a solution $U$ to \eqref{Eq1.5} satisfying
$\overline{u} \leq U \leq Z$.

\medskip
By Lemma \ref{l1.3}, $W := \Delta U > 0$, so $U$ is increasing. As $\Delta W = -U^{-1} < 0$, $W$ is decreasing and
$W(r)\geq Cr^2U^{-1}(r)$, see for example \eqref{1.2new}. By Corollary
\ref{C1.4}, $\lim_{r\to\infty}W(r)=\alpha>0.$ Therefore $\lim_{r\to\infty}\frac{U}{r^2}= \frac{\alpha}{2N}$ and
$$
\lim_{r\to\infty}rW'(r)=-\lim_{r\to\infty}r^{2-N} \int_{0}^{r}\frac{t^{N-1}}{U(t)}dt = -\lim_{r\to\infty}\frac{r^{2}}{(N - 2)U(r)}
 = -\frac{2N}{(N-2)\alpha} < 0.
$$
This implies that $W(r) < 0$ for $r$ large enough, which contradicts $W > 0$. \qed

\medskip
\noindent
{\bf Acknowledgements:} The authors are grateful to A. Farina for sending us the interesting preprint \cite{FF} and the papers \cite{Wa, Wa1} of Walter.

\bigskip
\noindent
Baishun Lai -
Institute of Contemporary Mathematics, Henan University, Kaifeng 475004, P.R. China.

\noindent
{\sl E-mail address}: {\tt laibaishun@henu.edu.cn}

\medskip
\noindent
Dong Ye - IECL, UMR 7502, D\'epartement de Math\'ematiques, Universit\'e de Lorraine, Ile de Saulcy, 57045 Metz, France.

\noindent
{\sl E-mail address}: {\tt dong.ye@univ-lorraine.fr}

\end{document}